\documentclass[11pt, a4paper]{article}
\usepackage{graphicx} 
\usepackage{geometry}
\usepackage{amsmath}
\usepackage{amssymb}
\usepackage{amsthm}
\usepackage{tikz}

\usepackage{authblk}
\usepackage{caption}
\usepackage{subcaption}

\title{Deep Neural Networks with Symplectic Preservation Properties}
\author[1]{Qing He}
\author[2]{Wei Cai}

\affil[1,2]{Dept. of Mathematics, Southern Methodist University, Dallas, TX, 75275}

\begin{document}
\maketitle

\begin{abstract}
We propose a deep neural network architecture designed such that its output forms an invertible symplectomorphism of the input. This design draws an analogy to the real-valued non-volume-preserving (real NVP) method used in normalizing flow techniques. Utilizing this neural network type allows for learning tasks on unknown Hamiltonian systems without breaking the inherent symplectic structure of the phase space.


\end{abstract}

\section{Introduction}
For an unknown Hamiltonian system, our objective is to learn the flow mapping over a fixed time period \( T \). Specifically, we seek to determine the map \( \Phi_T \) that computes \( (q, p)_{t=T} \) given an initial condition \( (q, p)_{t=0} = (q_0, p_0) \). Such problems arise, for instance, when analyzing a sequence of system snapshots at times \( 0, T, 2T, 3T, \ldots \). The key information we possess about this mapping is its property as a \textbf{symplectomorphism} (or \textbf{canonical transformation}), implying that the Jacobian of \( \Phi_T \) belongs to the symplectic group \( Sp(2n) \), where \( n \) is the dimensionality of the system's configuration space \cite{da2004lectures, goldstein1980classical}.

In this study, we propose a neural network structure designed to ensure that its output is precisely a symplectomorphism of the input. "Precisely" here means that the Jacobian of the mapping defined by the neural network is exactly a symplectic matrix, accounting only for minimal rounding errors inherent to floating-point arithmetic. Importantly, this framework eliminates the need to introduce an additional "deviation-from-symplecticity penalty term" in our learning objective because the inherent structure of the network guarantees that the symplectomorphism condition cannot be violated.

The approach draws inspiration from the \textit{real NVP} method \cite{Dinh2016DensityEU}, which is primarily used for density estimation of probability measures and differs significantly in purpose from our intended application. Nonetheless, this work leverages real NVP's elegant methodology for constructing explicitly invertible neural networks. The method we propose represents a "symplectic adaptation" of this technique, employing building blocks akin to those in real NVP while ensuring the preservation of symplecticity throughout. This adaptation involves replacing components that could potentially compromise the symplectic property of the mapping.

\section{Preliminaries}

\subsection{Symplectic Structures and Symplectomorphism}
On $\mathbb{R}^{2n}$, we denote the standard Cartesian coordinates as $q_{1} ,\cdots ,q_{n} ,p_{1} ,\cdots ,p_{n}$, corresponding to the "position" and "momentum" coordinates in Hamiltonian mechanics. The standard symplectic form on $\mathbb{R}^{2n}$ is the differential 2-form
\begin{equation}
\omega =\sum _{i=1}^{n}\mathrm{d} q_{i} \land \mathrm{d} p_{i} ,
\end{equation}
and a transformation $\varphi :\mathbb{R}^{2n}\rightarrow \mathbb{R}^{2n}$ is called a symplectomorphism if $\varphi ^{*} \omega =\omega $. This means
\begin{equation}
\sum _{i=1}^{n}\mathrm{d} Q_{i} \land \mathrm{d} P_{i} =\sum _{i=1}^{n}\mathrm{d} q_{i} \land \mathrm{d} p_{i} ,
\end{equation}
where 
\begin{equation}
( Q_{1} ,\cdots ,Q_{n} ,P_{1} ,\cdots ,P_{n}) =\varphi ( q_{1} ,\cdots ,q_{n} ,p_{1} ,\cdots ,p_{n}) ,
\end{equation}
or equivalently,
\begin{equation}
J_{\varphi }^{\top } \Omega J_{\varphi } =\Omega ,
\end{equation}
where 
\begin{equation}
J_{\varphi } =\begin{pmatrix}
\frac{\partial Q_{1}}{\partial q_{1}} & \cdots  & \frac{\partial Q_{1}}{\partial q_{n}} & \frac{\partial Q_{1}}{\partial p_{1}} & \cdots  & \frac{\partial Q_{1}}{\partial p_{n}}\\
\vdots  & \ddots  & \vdots  & \vdots  & \ddots  & \vdots \\
\frac{\partial Q_{n}}{\partial q_{1}} & \cdots  & \frac{\partial Q_{n}}{\partial q_{n}} & \frac{\partial Q_{n}}{\partial p_{1}} & \cdots  & \frac{\partial Q_{n}}{\partial p_{n}}\\
\frac{\partial P_{1}}{\partial q_{1}} & \cdots  & \frac{\partial P_{1}}{\partial q_{n}} & \frac{\partial P_{1}}{\partial p_{1}} & \cdots  & \frac{\partial P_{1}}{\partial p_{n}}\\
\vdots  & \ddots  & \vdots  & \vdots  & \ddots  & \vdots \\
\frac{\partial P_{n}}{\partial q_{1}} & \cdots  & \frac{\partial P_{n}}{\partial q_{n}} & \frac{\partial P_{n}}{\partial p_{1}} & \cdots  & \frac{\partial P_{n}}{\partial p_{n}}
\end{pmatrix}
\end{equation}
is the Jacobian matrix of $\varphi $, and 
\begin{equation}
\Omega =\begin{pmatrix}
0_{n\times n} & I_{n\times n}\\
-I_{n\times n} & 0_{n\times n}
\end{pmatrix}
\end{equation}
is the matrix of the standard symplectic form $\omega$.

The most essential property of a Hamiltonian system
\begin{equation} \label{HamiltSys}
\begin{cases}
\displaystyle \frac{\mathrm{d} q_{i}}{\mathrm{d} t} =\frac{\partial H}{\partial p_{i}}, & \\[10pt]
\displaystyle \frac{\mathrm{d} p_{i}}{\mathrm{d} t} =-\frac{\partial H}{\partial q_{i}}, & 
\end{cases} i=1,2,\cdots ,n,
\end{equation}
where
\begin{equation*}
    H = H(q_1, \cdots, q_n, p_1, \cdots, p_n, t) \in C^{2}(\mathbb{R}^{2n+1})
\end{equation*}
is that its flow map defines a family of symplectomorphisms. This means that if we solve (\ref{HamiltSys}) from time $t_0$ to time $t_1$, then the mapping defined by $(q(t_0), p(t_0)) \to (q(t_1), p(t_1))$ is an $\mathbb{R}^{2n} \to \mathbb{R}^{2n}$ symplectomorphism. The inverse is also true: If a differential equation system on $\mathbb{R}^{2n} \to \mathbb{R}^{2n}$ satisfies than the flow maps are symlectomorphisms, then there exists a function $H \in C^{2}(\mathbb{R}^{2n+1})$ such that the system can be written as Hamiltonian system (\ref{HamiltSys}). 

\subsubsection{Example: Shearing}

One simplest example of symplecticomorphism comes from the symplectic Euler method for separable Hamiltonian. Suppose $F:\mathbb{R}^{n}\rightarrow \mathbb{R}$ is a smooth function, then 
\begin{equation} \label{61XeKpBFWAGpC1mAJ52S}
\begin{cases}
Q_{i} =q_{i} & \\
P_{i} =p_{i} +\frac{\partial F}{\partial q_{i}}( q_{1} ,\cdots ,q_{n}) & 
\end{cases}
\end{equation}
is a symplectic transformation, because 
\begin{equation*}
\begin{aligned}
\sum _{i=1}^{n}\mathrm{d} Q_{i} \land \mathrm{d} P_{i} = & \sum _{i=1}^{n}\mathrm{d} q_{i} \land \mathrm{d}\left( p_{i} +\frac{\partial F}{\partial q_{i}}( q_{1} ,\cdots ,q_{n})\right)\\
= & \sum _{i=1}^{n}\mathrm{d} q_{i} \land \mathrm{d} p_{i} +\sum _{i=1}^{n}\mathrm{d} q_{i} \land \mathrm{d}\frac{\partial F}{\partial q_{i}}( q_{1} ,\cdots ,q_{n})\\
= & \sum _{i=1}^{n}\mathrm{d} q_{i} \land \mathrm{d} p_{i} -\mathrm{d}(\mathrm{d} F( q_{1} ,\cdots ,q_{n}))
\end{aligned}
\end{equation*}
and the result comes from the identity $\mathrm{d}(\mathrm{d} F) =0$. And similarly, 
\begin{equation} \label{etUqTh6SkoIVJZwDufc3}
\begin{cases}
Q_{i} =q_{i} +\frac{\partial G}{\partial p_{i}}( p_{1} ,\cdots ,p_{n}) & \\
P_{i} =p_{i} & 
\end{cases}
\end{equation}
is also a symplectomorphism, where $G:\mathbb{R}^{n}\rightarrow \mathbb{R}$ is a smooth function. We call the symplectomorphism given by (\ref{61XeKpBFWAGpC1mAJ52S}) or (\ref{etUqTh6SkoIVJZwDufc3}) a \textbf{symplectic shearing}.

\subsubsection{Example: Stretching}

Another example is the "coordinate stretching" transformation. A diagonal linear transformation on $\mathbb{R}^{2n}$ is symplectic if and only if it has the form 
\begin{equation} \label{L3SLMeGtQT0E2b5kwMMF}
( q_{1} ,\cdots ,q_{n} ,p_{1} ,\cdots ,p_{n}) \mapsto \left( k_{1} q_{1} ,\cdots ,k_{n} q_{n} ,\frac{p_{1}}{k_{1}} ,\cdots ,\frac{p_{n}}{k_{n}}\right) ,
\end{equation}
where $k_{1} ,\cdots ,k_{n}$ are nonzero constants. Now we make it more general, supposing that each $k_{i}$'s are functions of the coordinates $q_{1} ,\cdots ,q_{n} ,p_{1} ,\cdots ,p_{n}$. Then
\begin{equation} \label{V7y4EJWAWgqp6kQEVCPW}
\begin{aligned}
\sum _{i=1}^{n}\mathrm{d}( k_{i} q_{i}) \land \mathrm{d}\frac{p_{i}}{k_{i}} = & \sum _{i=1}^{n}( k_{i}\mathrm{d} q_{i} +q_{i}\mathrm{d} k_{i}) \land \left(\frac{\mathrm{d} p_{i}}{k_{i}} -\frac{p_{i}\mathrm{d} k_{i}}{k_{i}^{2}}\right)\\
= & \sum _{i=1}^{n}\mathrm{d} q_{i} \land \mathrm{d} p_{i} +\frac{q_{i}}{k_{i}}\mathrm{d} k_{i} \land \mathrm{d} p_{i} -\frac{p_{i}\mathrm{d} q_{i} \land \mathrm{d} k_{i}}{k_{i}} +0\\
= & \sum _{i=1}^{n}\mathrm{d} q_{i} \land \mathrm{d} p_{i} -\frac{q_{i}\mathrm{d} p_{i} +p_{i}\mathrm{d} q_{i}}{k_{i}} \land \mathrm{d} k_{i}\\
= & \sum _{i=1}^{n}\mathrm{d} q_{i} \land \mathrm{d} p_{i} -\frac{\mathrm{d}( p_{i} q_{i})}{k_{i}} \land \mathrm{d} k_{i} ,
\end{aligned}
\end{equation}
therefore, a transformation given as (\ref{L3SLMeGtQT0E2b5kwMMF}) is symplectic if and only if the condition
\begin{equation} \label{hWrVPnxWXeBR4lBhTxfE}
\sum _{i=1}^{n}\frac{\mathrm{d}( p_{i} q_{i})}{k_{i}} \land \mathrm{d} k_{i} =0
\end{equation}
is satisfied, the mapping (\ref{L3SLMeGtQT0E2b5kwMMF}) is symplectic. Note that (\ref{hWrVPnxWXeBR4lBhTxfE}) can be written as 
\begin{equation*}
\sum _{i=1}^{n}\mathrm{d}( p_{i} q_{i}) \land \mathrm{d}\ln |k_{i} |=0,
\end{equation*}
and accoring to Poincaré's Lemma, (\ref{hWrVPnxWXeBR4lBhTxfE}) is satisfied if 
\begin{equation} \label{kTInCPyvuMTkv1XGplqK}
\sum _{i=1}^{n}\ln |k_{i} |\mathrm{d}( p_{i} q_{i}) =\mathrm{d} \varphi 
\end{equation}
for some smooth function $\varphi :\mathbb{R}^{2n}\rightarrow \mathbb{R}$. The condition (\ref{kTInCPyvuMTkv1XGplqK}) is satisfied when $\varphi $ can be expressed as 
\begin{equation*}
\varphi ( q_{1} ,\cdots ,q_{n} ,p_{1} ,\cdots ,p_{n}) =\Phi ( p_{1} q_{1} ,p_{2} q_{2} ,\cdots ,p_{n} q_{n})
\end{equation*}
for some $\Phi :\mathbb{R}^{n}\rightarrow \mathbb{R}$, and 
\begin{equation} \label{gBrP1tE0YYDAIsdyZWlJ}
k_{i} =\pm \mathrm{e}^{\Phi _{i}( p_{1} q_{1} ,p_{2} q_{2} ,\cdots ,p_{n} q_{n})}
\end{equation}
holds, where $\Phi _{i}$ is the partial derivative of $\Phi $ on its $i$-ith argument:
\begin{equation} \label{eKXfv2VqyDx49nbPufw5}
\Phi _{i}( x_{1} ,\cdots ,x_{n}) =\frac{\partial \Phi }{\partial x_{i}}( x_{1} ,\cdots ,x_{n}) .
\end{equation}
We call the symplectomorihism given by (\ref{L3SLMeGtQT0E2b5kwMMF}) and (\ref{gBrP1tE0YYDAIsdyZWlJ}) a \textbf{symplectic stretching}.

\subsection{Real NVP}

Real NVP (Real-valued Non-Volume Preserving) \cite{Dinh2016DensityEU, bishop2023deep} is a generative model used for density estimation. Real NVP networks use invertible transformations, allowing us to go back and forth between the original and transformed spaces. The structure of real NVP is as follows: The input and output of the network are both $N$-dimensional vectors. An $N$-dimensional vector
\begin{equation*}
z=( z_{1} ,z_{2} ,\cdots ,z_{N})
\end{equation*}
received as the input is partitioned in to two parts 
\begin{equation*}
z=(\underbrace{z_{1} ,\cdots ,z_{n}}_{A} ,\underbrace{z_{n+1} ,\cdots ,z_{N}}_{B}) :=( z_{A} ,z_{B}) .
\end{equation*}
A Real NVP transformation keeps one of the parts unchanged and perform an "entry-wise linear transformation" on the other part, whose coefficients are determined by the unchanged part. Specifically, the input $z$ undergoes the following transformation:
\begin{equation} \label{NlCG4CGPyQTqtKx5RvwW}
\begin{cases}
x_{A} =z_{A} & \\
x_{B} =\mathrm{e}^{s( z_{A})} \odot z_{B} +b( z_{A}) & 
\end{cases}
\end{equation}
where $s,b:\mathbb{R}^{n}\rightarrow \mathbb{R}^{N-n}$ are two functions which are given as a neural networks in practice, and the symbol "$\odot $" the Hadamard product (entry-wise product) operator:
\begin{equation*}
( x_{1} ,\cdots ,x_{n}) \odot ( y_{1} ,\cdots ,y_{n}) =( x_{1} y_{1} ,\cdots ,x_{n} y_{n}) .
\end{equation*}
The inverse of this mapping (\ref{NlCG4CGPyQTqtKx5RvwW}) is clear:
\begin{equation} \label{HLYtePSjTjYT1qEkNvpI}
\begin{cases}
z_{A} =x_{A} & \\
z_{B} =\mathrm{e}^{-s( z_{A})} \odot ( x_{B} -b( x_{A})) . & 
\end{cases}
\end{equation}
The transformation (\ref{NlCG4CGPyQTqtKx5RvwW}) is often exhibited as a diagram like .

\begin{figure}[htbp]
    \centering
    \includegraphics[height=0.25\textwidth]{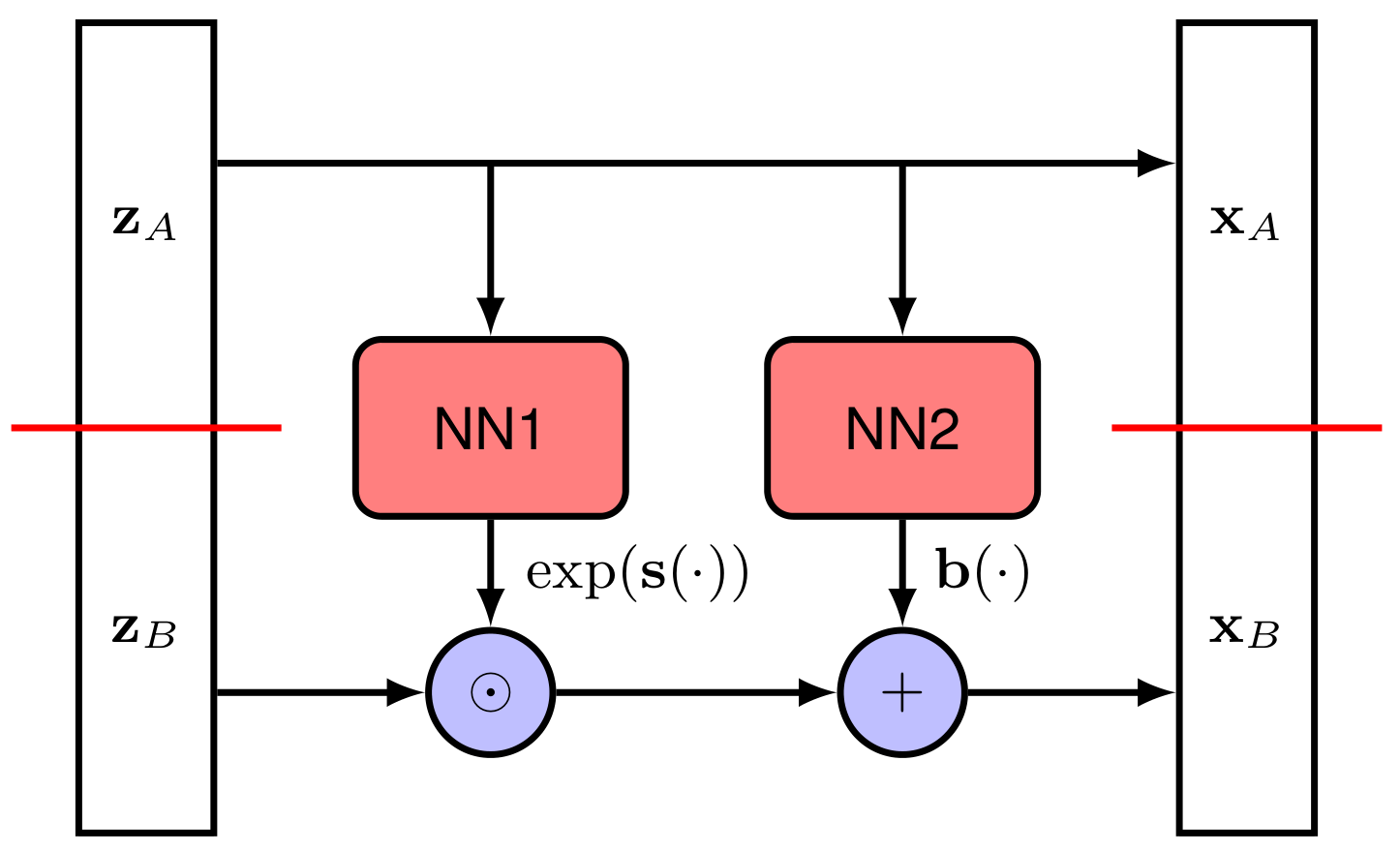}
    \caption{A diagram of the transformation (\ref{NlCG4CGPyQTqtKx5RvwW})}
    \label{fig:real_nvp1}
\end{figure}

The apparent limitation of transformation (\ref{NlCG4CGPyQTqtKx5RvwW}) is that it does not change the part $z_{A}$. This can be quickly fixed by appending another real NVP block that keeps the $x_{B}$ part unchanged: 
\begin{equation} \label{Y9lpGfmCh2alj0m1WYpk}
\begin{cases}
y_{A}\leftarrow \mathrm{e}^{\tilde{s}( x_{A})} \odot x_{B} +\tilde{b}( x_{A}) & \\
y_{B}\leftarrow x_{B} & 
\end{cases}
\end{equation}
where $\tilde{s} ,\tilde{b} :\mathbb{R}^{N-n}\rightarrow \mathbb{R}^{n}$ are another two neural network functions, so the composed transformation from $z$ to $y$ given by (\ref{NlCG4CGPyQTqtKx5RvwW}) and (\ref{Y9lpGfmCh2alj0m1WYpk}) do not keep any component unchanged. This can be exhibited as a diagram like .

\begin{figure}[htbp]
    \centering
    \includegraphics[height=0.25\textwidth]{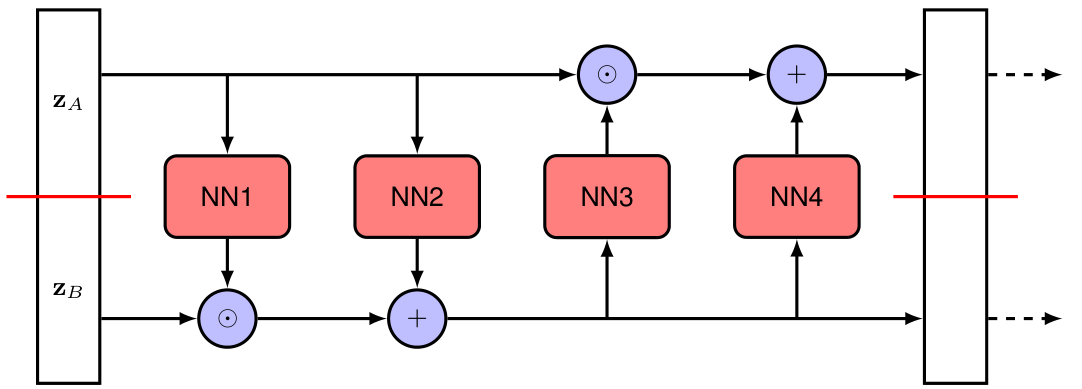}
    \caption{A diagram of the composition transformation}
    \label{fig:real_nvp2}
\end{figure}

Of course, we can stack more layers like this to improve the expressivity of the network.

\section{Symplectomorphism Neural Network (SymplectoNet, SpNN)}
\subsection{Structure}

For our goal of building symplectomorphism neural network, the problem of real NVP is directly exhibited in its name: "NVP" means "non-volume-preserving", while a symplectomorphism has to be volume preserving. Indeed, to make real NVP volume preserving (from "real NVP" to "real VP"), there is a quick fixation: one only needs to add an extra layer
\begin{equation*}
( s_{1} ,\cdots ,s_{N})\rightarrow ( s_{1} ,\cdots ,s_{N}) -\overline{s}( 1,\cdots ,1) ,\quad \overline{s} =\frac{1}{N}\sum _{i=1}^{N} s_{i}
\end{equation*}
after the output layer of the network that subtracts the average of the network. Unfortunately, mere volume-preserving property does not guarantee symplecticity. We need further adjustments.

Indeed, we can decompose (\ref{NlCG4CGPyQTqtKx5RvwW}) into two transformations: a "stretching"
\begin{equation} \label{TxKO6T1S5Oh6u1d74FaP}
\begin{cases}
\xi _{A} =z_{A} & \\
\xi _{B} =\mathrm{e}^{s( z_{A})} \odot z_{B} , & 
\end{cases}
\end{equation}
and a "shearing"
\begin{equation} \label{k7u8zYsKbYW5JgjaJEWe}
\begin{cases}
x_{A} =\xi _{A} & \\
x_{B} =\xi _{B} +b( \xi _{A}) . & 
\end{cases}
\end{equation}
Neither of these two transformations are guaranteed to be symplectic. Nevertheless, we have introduced their symplectic counterparts in the last section: Indeed, we can write (\ref{61XeKpBFWAGpC1mAJ52S}), (\ref{etUqTh6SkoIVJZwDufc3}) and (\ref{L3SLMeGtQT0E2b5kwMMF}) (where (\ref{gBrP1tE0YYDAIsdyZWlJ}), (\ref{eKXfv2VqyDx49nbPufw5})) into a more compact form 
\begin{equation} \label{crgEDid00Phhk3uChgTs}
\begin{cases}
Q=q & \\
P=p+\nabla F( q) , & 
\end{cases}
\end{equation}
\begin{equation} \label{oprjyQo8NkwkcJbVWpDG}
\begin{cases}
Q=q+\nabla G( p) & \\
P=p, & 
\end{cases}
\end{equation}
\begin{equation} \label{e2tglIQQmVwg3ghbgXp5}
\begin{cases}
Q=\mathrm{e}^{\nabla \Phi ( q\odot p)} \odot q & \\
P=\mathrm{e}^{-\nabla \Phi ( q\odot p)} \odot p, & 
\end{cases}
\end{equation}
where $q=( q_{1} ,\cdots ,q_{n})$, $p=( p_{1} ,\cdots ,p_{n})$, $Q=( Q_{1} ,\cdots ,Q_{n})$, $P=( P_{1} ,\cdots ,P_{n})$. And "$\odot $" is the Hadamard product as before. And now their correspondence with (\ref{TxKO6T1S5Oh6u1d74FaP}) and (\ref{k7u8zYsKbYW5JgjaJEWe}) are clear: (\ref{e2tglIQQmVwg3ghbgXp5}) is exactly (\ref{k7u8zYsKbYW5JgjaJEWe}) when $\dim x_{A}$ and $\dim x_{B}$ are of the same dimension, and $b$ is the gradient of a function; while (\ref{e2tglIQQmVwg3ghbgXp5}) is a symmetrized version of (\ref{TxKO6T1S5Oh6u1d74FaP}):
\begin{equation*}
\begin{cases}
\xi _{A} =\mathrm{e}^{-s( z_{A} \odot z_{B})} \odot z_{A} & \\
\xi _{B} =\mathrm{e}^{s( z_{A} \odot z_{B})} \odot z_{B} , & 
\end{cases}
\end{equation*}
with $s$ being the gradient of a function. We denote the transformations defined by (\ref{crgEDid00Phhk3uChgTs}), (\ref{oprjyQo8NkwkcJbVWpDG}), (\ref{e2tglIQQmVwg3ghbgXp5}) as $\operatorname{qSh}_{F}$, $\operatorname{pSh}_{G}$ and $\operatorname{St_{\Phi }}$, which are short hands for "q-shearing", "p-shearing" and "stretching", respectively. These becomes the basic building blocks of the "symplectic version of real NVP" once we take $F$, $G$ and $\Phi $ in these transformations as trainable neural networks.

Now we have introduced all the basic symplecticomorphism building blocks, and a \textbf{symplectomorphism neural network} (\textbf{SymplectoNet}, or even shorter, \textbf{SpNN}) is \textbf{a neural network designed as an arbitrary finite composition of }$\operatorname{qSh}_{F}$\textbf{, }$\operatorname{pSh}_{G}$\textbf{ and }$\operatorname{St_{\Phi }}$\textbf{ where }$F$\textbf{, }$G$\textbf{ and }$\Phi $\textbf{ are arbitrary neural networks with }$n$\textbf{-dimensional input and one-dimensional output}.

\begin{figure}[htbp]
    \centering
    \tikzset{every picture/.style={line width=0.75pt}} 

\begin{tikzpicture}[x=0.75pt,y=0.75pt,yscale=-1,xscale=1]

\draw   (102,1.37) -- (129.84,1.37) -- (129.84,101.43) -- (102,101.43) -- cycle ;
\draw   (102,101.43) -- (129.84,101.43) -- (129.84,201.5) -- (102,201.5) -- cycle ;
\draw    (129.9,51.15) -- (376.04,51.15) ;
\draw    (130.18,151.16) -- (178.48,151.16) ;
\draw    (189.76,51.4) -- (189.76,82.48) ;
\draw  [fill={rgb, 255:red, 208; green, 2; blue, 27 }  ,fill opacity=0.4 ] (172.51,82.55) -- (206.9,82.55) -- (206.9,116.94) -- (172.51,116.94) -- cycle ;
\draw    (189.76,116.73) -- (189.76,139.8) ;
\draw  [fill={rgb, 255:red, 189; green, 16; blue, 224 }  ,fill opacity=0.4 ] (178.48,151.39) .. controls (178.48,145.16) and (183.53,140.11) .. (189.76,140.11) .. controls (195.99,140.11) and (201.04,145.16) .. (201.04,151.39) .. controls (201.04,157.62) and (195.99,162.67) .. (189.76,162.67) .. controls (183.53,162.67) and (178.48,157.62) .. (178.48,151.39) -- cycle ;
\draw    (201.04,151.14) -- (376.04,151.14) ;
\draw    (184.6,151.2) -- (195.27,151.2) ;
\draw    (189.93,156.7) -- (189.93,146.03) ;
\draw    (268.93,51.15) -- (268.93,89.73) ;
\draw    (268.93,112.57) -- (268.93,151.14) ;
\draw  [fill={rgb, 255:red, 189; green, 16; blue, 224 }  ,fill opacity=0.4 ] (257.65,101.01) .. controls (257.65,94.78) and (262.7,89.73) .. (268.93,89.73) .. controls (275.16,89.73) and (280.21,94.78) .. (280.21,101.01) .. controls (280.21,107.24) and (275.16,112.29) .. (268.93,112.29) .. controls (262.7,112.29) and (257.65,107.24) .. (257.65,101.01) -- cycle ;
\draw   (262.81,101.01) .. controls (262.81,97.62) and (265.55,94.88) .. (268.93,94.88) .. controls (272.31,94.88) and (275.06,97.62) .. (275.06,101.01) .. controls (275.06,104.39) and (272.31,107.13) .. (268.93,107.13) .. controls (265.55,107.13) and (262.81,104.39) .. (262.81,101.01) -- cycle ;
\draw   (269.43,101.01) .. controls (269.43,100.73) and (269.21,100.51) .. (268.93,100.51) .. controls (268.65,100.51) and (268.43,100.73) .. (268.43,101.01) .. controls (268.43,101.28) and (268.65,101.51) .. (268.93,101.51) .. controls (269.21,101.51) and (269.43,101.28) .. (269.43,101.01) -- cycle ;
\draw  [fill={rgb, 255:red, 208; green, 2; blue, 27 }  ,fill opacity=0.4 ] (290.47,82.94) -- (324.86,82.94) -- (324.86,117.33) -- (290.47,117.33) -- cycle ;
\draw    (280.21,101.01) -- (289.97,101.01) ;
\draw    (325.18,101.24) -- (336.26,101.24) ;
\draw  [fill={rgb, 255:red, 208; green, 2; blue, 27 }  ,fill opacity=0.4 ] (336.1,90.55) -- (364.55,90.55) -- (364.55,112.18) -- (336.1,112.18) -- cycle ;
\draw  [fill={rgb, 255:red, 189; green, 16; blue, 224 }  ,fill opacity=0.4 ] (376.04,151.22) .. controls (376.04,144.99) and (381.09,139.94) .. (387.32,139.94) .. controls (393.55,139.94) and (398.6,144.99) .. (398.6,151.22) .. controls (398.6,157.45) and (393.55,162.5) .. (387.32,162.5) .. controls (381.09,162.5) and (376.04,157.45) .. (376.04,151.22) -- cycle ;
\draw   (381.2,151.22) .. controls (381.2,147.84) and (383.94,145.1) .. (387.32,145.1) .. controls (390.71,145.1) and (393.45,147.84) .. (393.45,151.22) .. controls (393.45,154.6) and (390.71,157.35) .. (387.32,157.35) .. controls (383.94,157.35) and (381.2,154.6) .. (381.2,151.22) -- cycle ;
\draw   (387.82,151.22) .. controls (387.82,150.94) and (387.6,150.72) .. (387.32,150.72) .. controls (387.05,150.72) and (386.82,150.94) .. (386.82,151.22) .. controls (386.82,151.5) and (387.05,151.72) .. (387.32,151.72) .. controls (387.6,151.72) and (387.82,151.5) .. (387.82,151.22) -- cycle ;
\draw  [fill={rgb, 255:red, 189; green, 16; blue, 224 }  ,fill opacity=0.4 ] (376.04,50.93) .. controls (376.04,44.7) and (381.09,39.65) .. (387.32,39.65) .. controls (393.55,39.65) and (398.6,44.7) .. (398.6,50.93) .. controls (398.6,57.16) and (393.55,62.22) .. (387.32,62.22) .. controls (381.09,62.22) and (376.04,57.16) .. (376.04,50.93) -- cycle ;
\draw   (381.2,50.93) .. controls (381.2,47.55) and (383.94,44.81) .. (387.32,44.81) .. controls (390.71,44.81) and (393.45,47.55) .. (393.45,50.93) .. controls (393.45,54.32) and (390.71,57.06) .. (387.32,57.06) .. controls (383.94,57.06) and (381.2,54.32) .. (381.2,50.93) -- cycle ;
\draw   (387.82,50.93) .. controls (387.82,50.66) and (387.6,50.43) .. (387.32,50.43) .. controls (387.05,50.43) and (386.82,50.66) .. (386.82,50.93) .. controls (386.82,51.21) and (387.05,51.43) .. (387.32,51.43) .. controls (387.6,51.43) and (387.82,51.21) .. (387.82,50.93) -- cycle ;
\draw    (387.32,62.22) -- (387.28,71.71) ;
\draw    (387.32,94.17) -- (387.32,139.94) ;
\draw    (364.93,101.24) -- (387.05,101.24) ;
\draw  [fill={rgb, 255:red, 189; green, 16; blue, 224 }  ,fill opacity=0.4 ] (373.85,72.17) -- (400.8,72.17) -- (400.8,93.83) -- (373.85,93.83) -- cycle ;
\draw    (398.6,50.93) -- (459.6,50.93) ;
\draw    (398.6,151.22) -- (529.8,151.22) ;
\draw    (471.26,151.17) -- (471.26,120.1) ;
\draw  [fill={rgb, 255:red, 208; green, 2; blue, 27 }  ,fill opacity=0.4 ] (454.01,120.03) -- (488.4,120.03) -- (488.4,85.64) -- (454.01,85.64) -- cycle ;
\draw    (471.26,85.85) -- (471.26,62.77) ;
\draw  [fill={rgb, 255:red, 189; green, 16; blue, 224 }  ,fill opacity=0.4 ] (459.98,51.18) .. controls (459.98,57.41) and (465.03,62.46) .. (471.26,62.46) .. controls (477.49,62.46) and (482.54,57.41) .. (482.54,51.18) .. controls (482.54,44.95) and (477.49,39.9) .. (471.26,39.9) .. controls (465.03,39.9) and (459.98,44.95) .. (459.98,51.18) -- cycle ;
\draw    (466.1,51.37) -- (476.77,51.37) ;
\draw    (471.43,45.87) -- (471.43,56.54) ;
\draw    (482.54,51.18) -- (529.3,51.18) ;
\draw   (529.5,1.37) -- (557.34,1.37) -- (557.34,101.43) -- (529.5,101.43) -- cycle ;
\draw   (529.5,101.43) -- (557.34,101.43) -- (557.34,201.5) -- (529.5,201.5) -- cycle ;
\draw  [dash pattern={on 4.5pt off 4.5pt}] (154.76,31.4) -- (224.76,31.4) -- (224.76,181.85) -- (154.76,181.85) -- cycle ;
\draw  [dash pattern={on 4.5pt off 4.5pt}] (248.26,31.4) -- (416.3,31.4) -- (416.3,181.85) -- (248.26,181.85) -- cycle ;
\draw  [dash pattern={on 4.5pt off 4.5pt}] (435.76,30.9) -- (505.76,30.9) -- (505.76,181.35) -- (435.76,181.35) -- cycle ;

\draw (115.92,51.4) node    {$q$};
\draw (115.92,151.47) node    {$p$};
\draw (189.7,99.74) node    {$\nabla F$};
\draw (307.67,100.13) node    {$\nabla \Phi $};
\draw (350.33,101.36) node   [align=left] {exp};
\draw (387.33,83) node   [align=left] {$1/x$};
\draw (471.2,102.83) node    {$\nabla G$};
\draw (191.18,23.2) node [anchor=south] [inner sep=0.75pt]    {$\operatorname{qSh}_{F}$};
\draw (472.18,22.8) node [anchor=south] [inner sep=0.75pt]    {$\operatorname{pSh}_{G}$};
\draw (330.18,23.2) node [anchor=south] [inner sep=0.75pt]    {$\operatorname{St}_{\Phi }$};

\end{tikzpicture}
    \caption{The diagram expression of $\operatorname{pSh}_{G} \circ \operatorname{St}_{\Phi } \circ \operatorname{qSh}_{F}$}
    \label{fig:spnnstruct}
\end{figure}
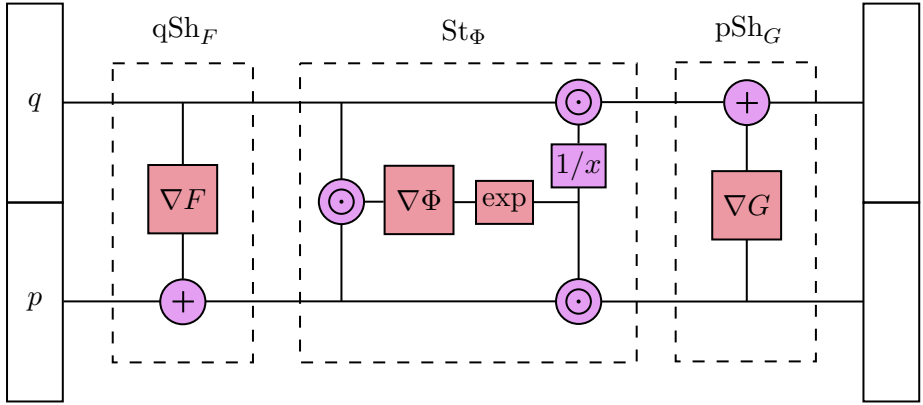

Of course, the expressivity of this network depends on the complexity of the underlying neural networks $F$, $G$ and $H$, and also on the number of the building blocks we stacked. Indeed, the latter can be even more essential: e.g. if we only use less than four symplectic shearing blocks, we cannot even cover all the linear symplectomorphisms no matter how complicated the underlying network $F$ and $G$ are, because the Jacobian of a shearing transformation is of the form 
\begin{equation*}
\begin{pmatrix}
I & \\
B & I
\end{pmatrix} \quad \text{or} \quad \begin{pmatrix}
I & C\\
 & I
\end{pmatrix} ,
\end{equation*}
where $B$, $C$ are symmetric $n\times n$ matrices. The degree of freedom of these matrices are $n( n+1) /2$, while $\dim Sp( 2n) =n( 2n+1)$, which is greater than $3n( n+1) /2$ for $n >1$. This is why I also designed the symplectic stretching layer $\operatorname{Str}_{\Phi }$. A good practice is to include both the p, q-shearing and the symplectic stretching layers in the network for at lease once. A simplest example is a network with structure $\operatorname{pSh}_{G} \circ \operatorname{St}_{\Phi } \circ \operatorname{qSh}_{F}$ (see \ref{fig:spnnstruct}), which is similar to the structure of a real NVP.

\subsection{SymplectoNet as Invertible Neural Network (INN)}

One of the most important features of real NVP is that it is explicitly invertible: one can write out (or, in a more techical term, build the computation graph of) the explicit expression of the neural network function's inverse function \cite{ishikawa2023universal}. Our SymplectoNet is inspired by real NVP, so a natural question is whether the SymplectoNet structure is explicitly invertible like real NVP. Next, we will show that the answer is yes.

Indeed, since the inverse of a composed function $f_{1} \circ f_{2} \circ \cdots \circ f_{k}$ is $f_{k}^{-1} \circ \cdots \circ f_{2}^{-1} \circ f_{1}^{-1}$, so we only need to prove that the basic building blocks, $\operatorname{pSh}_{G}$, $\operatorname{qSh}_{F}$ and $\operatorname{St}_{\Phi }$ are explicitly invertible. The inverse of $\operatorname{pSh}_{G}$, $\operatorname{qSh}_{F}$ are obvious: (\ref{crgEDid00Phhk3uChgTs}) is equivalent to 
\begin{equation*}
\begin{cases}
q=Q & \\
p=P-\nabla F( Q) , & 
\end{cases}
\end{equation*}
(\ref{oprjyQo8NkwkcJbVWpDG}) is equivalent to
\begin{equation*}
\begin{cases}
q=Q-\nabla G( P) & \\
p=P, & 
\end{cases}
\end{equation*}
therefore the inverse of $\operatorname{pSh}_{G}$, $\operatorname{qSh}_{F}$ are \ $\operatorname{pSh}_{-G}$, $\operatorname{qSh}_{-F}$, respectively. And finally we look at $\operatorname{St}_{\Phi }$. Notice that from (\ref{e2tglIQQmVwg3ghbgXp5}), we have 
\begin{equation*}
Q\odot P=\mathrm{e}^{\nabla \Phi ( q\odot p)} \odot q\odot \mathrm{e}^{-\nabla \Phi ( q\odot p)} \odot p=q\odot p,
\end{equation*}
 therefore 
\begin{equation} \label{62LLv9yMxX6WDegCRhk3}
\begin{cases}
q=\mathrm{e}^{-\nabla \Phi ( q\odot p)} \odot Q=\mathrm{e}^{-\nabla \Phi ( Q\odot P)} \odot Q & \\
p=\mathrm{e}^{\nabla \Phi ( q\odot p)} \odot P=\mathrm{e}^{\nabla \Phi ( Q\odot P)} \odot P, & 
\end{cases}
\end{equation}
this shows that the inverse of $\operatorname{St}_{\Phi }$ is exactly $\operatorname{St}_{-\Phi }$. In conclusion, we have 
\begin{equation} \label{jAml09tFVhozwfysrXRU}
\begin{cases}
\left(\operatorname{pSh}_{G}\right)^{-1} =\ \operatorname{pSh}_{-G} , & \\
\left(\operatorname{qSh}_{F}\right)^{-1} =\operatorname{qSh}_{-F} , & \\
\left(\operatorname{St}_{\Phi }\right)^{-1} =\operatorname{St}_{-\Phi } . & 
\end{cases}
\end{equation}
These results give a neat expression of inverting the SymplectoNet. E.g. the inverse of the SymplectoNet 
\begin{equation} \label{H4y2FzGoleGFTMfbR8lw}
\left(\operatorname{pSh}_{G} \circ \operatorname{St}_{\Phi } \circ \operatorname{qSh}_{F}\right)^{-1} =\operatorname{qSh}_{-F} \circ \operatorname{St}_{-\Phi } \circ \operatorname{pSh}_{-G} .
\end{equation}
This shows that the inverse of SymplectoNet is explicitly available.

\section{Extension to Family of Symplectomorphism}
A natural extension of the symplectomorphism neural network is to include some parameters $\tau_1, \tau_2, \cdots, \tau_K$ other that the canonical variables as inputs. This is can be easily achieved by changing the $F(q), G(p), \Phi(z)$ in the basic building blocks $\operatorname{qSh}_{F}$, $\operatorname{pSh}_{G}$ and $\operatorname{St}_{\Phi }$ into ($n + K$)-variable functions $F(q; \tau)$, $G(p; \tau)$, $\Phi(z; \tau)$, where $\tau = (\tau_1, \cdots, \tau_K)$, and modify the blocks given by (\ref{crgEDid00Phhk3uChgTs}) \textasciitilde (\ref{e2tglIQQmVwg3ghbgXp5}) into 
\begin{equation} \label{crgE34fdfk3uChgTs}
\begin{cases}
Q=q & \\
P=p+\nabla_q F( q, \tau) , & 
\end{cases}
\end{equation}
\begin{equation} \label{oprjdfgfdkcJbVWpDG}
\begin{cases}
Q=q+\nabla_p G( p, \tau) & \\
P=p, & 
\end{cases}
\end{equation}
\begin{equation} \label{e2sdfdfd3ghbgXp5}
\begin{cases}
Q=\mathrm{e}^{\nabla_z \Phi ( q\odot p, \tau)} \odot q & \\
P=\mathrm{e}^{-\nabla_z \Phi ( q\odot p, \tau)} \odot p, & 
\end{cases}
\end{equation}
With this modification, the network receives ($2n+K$)-dimensional vectors 
$$
(q_1,\cdots, q_n, p_1, \cdots, p_n, \tau_1, \cdots, \tau_K)
$$
as inputs and the output dimension is still $2n$, and for each fixed $\tau_1, \cdots, \tau_K$, the output vector is a symplectomorphism of the canonical part of the input vector, i.e. $(q_1,\cdots, q_n, p_1, \cdots, p_n)$. Thus, each choice of the parameters $\tau_1, \cdots, \tau_K$ defines a symplectomorphism, or we can say that the network defines a continuous family of symplectomorphisms parameterized by $\tau_1, \cdots, \tau_K$. A particularly common situation of this is when $K = 1$ and $\tau_1 = t$ represents the time variable. In this case, the network function can represent the solution of some Hamiltonian equation, and thanks to the symplectic property, of the network, there exists a Hamiltonian function 
$$
H = H(q_1, \cdots, q_n, p_1, \cdots, p_n, t)
$$
such that the network function represents \textit{exactly} the solution of its corresponding Hamiltonian system (\ref{HamiltSys}). Nevertheless, it is not guaranteed that the symplectomorphism family parameterized by $t$ forms a single-parameter symplectomorphism group, i.e. the corresponding Hamiltonian $H$ has to depend explicitly on time, and we do not have method to exactly cancel this dependency. 

By including more parameters (i.e. $K > 1$), it is also possible to apply this network for optimal control problems involving Hamiltonian dynamics.

\section{Some Preliminary Results}
\subsection{A Polar Nonlinear Mapping}

This example is learning a symplectic map
\begin{equation} \label{rvSHvBdfv8WquGXG0ixV}
( q,\ p)\rightarrow \left(\sqrt{2q}\cos p,\sqrt{2q}\sin p\right) =:( Q,P)
\end{equation}
A network with structure
\begin{equation*}
\operatorname{qSh}_{F_{1}} \circ \operatorname{pSh}_{G_{1}} \circ \operatorname{St}_{\Phi } \circ \operatorname{qSh}_{F_{2}} \circ \operatorname{pSh}_{G_{2}} ,
\end{equation*}
where $F_{1} ,G_{1} ,F_{2} ,G_{2}$ are $( 2,20,10,1)$ dense neural networks, and $\Phi $ is $( 2,10,1)$ dense neural network. The loss is the ordinary MSE loss. Adamax with learning rate 0.25 is applied here, and decay by factor 0.99 every 100 epoch. 

\begin{figure}
    \centering
    \begin{subfigure}{0.45\textwidth}
        \centering
        \includegraphics[height=0.2\textheight]{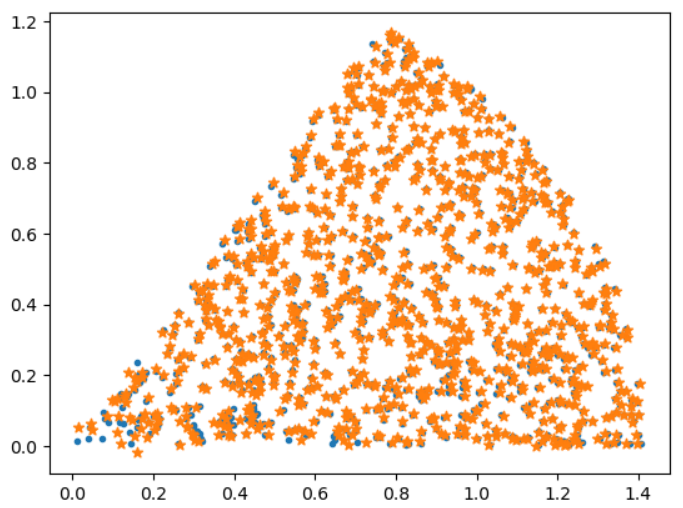}
        \caption{}
        \label{fig:0101}
    \end{subfigure}
    \begin{subfigure}{0.45\textwidth}
        \centering
        \includegraphics[height=0.2\textheight]{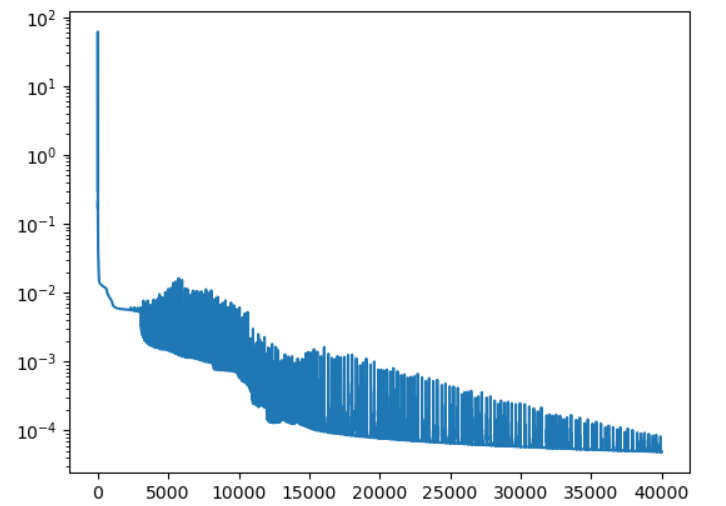}
        \caption{}
        \label{fig:0101_loss}
    \end{subfigure}
    \begin{subfigure}{0.45\textwidth}
        \centering
        \includegraphics[height=0.2\textheight]{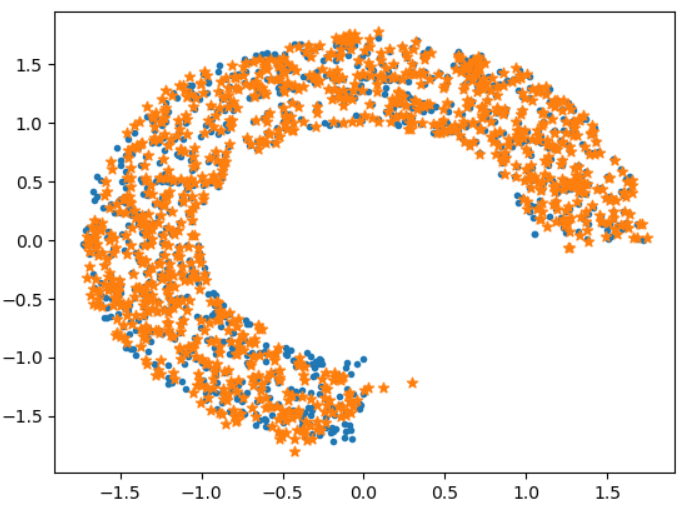}
        \caption{}
        \label{fig:123203pi2}
    \end{subfigure}
    \begin{subfigure}{0.45\textwidth}
        \centering
        \includegraphics[height=0.2\textheight]{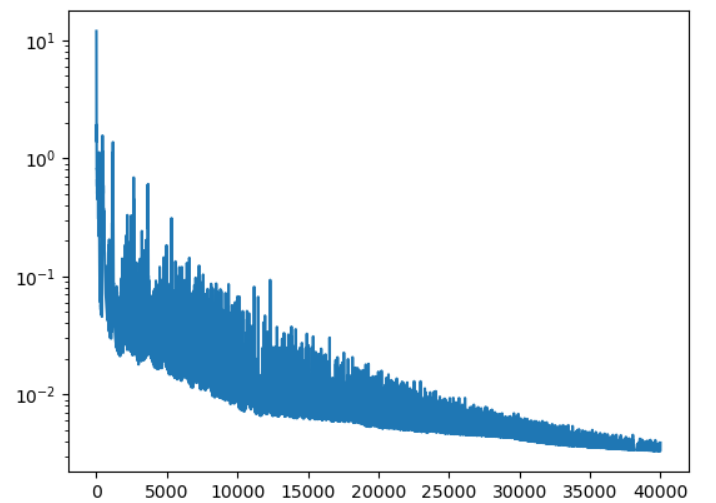}
        \caption{}
        \label{fig:123203pi2_loss}
    \end{subfigure}
    \caption{Numerical experiment results of symplectomorphism neural network fitting the symplectomorphism (\ref{rvSHvBdfv8WquGXG0ixV}). (a): The result of (\ref{rvSHvBdfv8WquGXG0ixV}) with $( q,p) \in [ 0,1] \times [ 0,1]$. Blue dots: true data; Orange stars: predicted results. Note that most of the error comes from data near $q=0$ because there is a singularity there; (b): The loss decay of (a); (c): The result of (\ref{rvSHvBdfv8WquGXG0ixV}) with $( q,p) \in [ 1/2,3/2] \times [ 0,3\pi / 2]$. Blue dots: true data; Orange stars: predicted results. Note that most of the error comes from data near $q=0$ because there is a singularity there; (d): The loss decay of (c);}
    \label{fig:spnn_for_polar_map}
\end{figure}

Firstly, some uniformly random points for 
\begin{equation*}
( q,p) \in [ 0,1] \times [ 0,1]
\end{equation*}
is sampled. The training went for 40,000 epochs, and the loss dropped from $0.3$ to about $10^{-5}$, and the plot is shown in Figure \ref{fig:0101}, and the loss decay is shown in Figure \ref{fig:0101_loss}.

Anoter numerical experiments concerning also (\ref{rvSHvBdfv8WquGXG0ixV}) but the domain changed to
\begin{equation*}
( q,p) \in \left[\frac{1}{2} ,\frac{3}{2}\right] \times \left[ 0,\frac{3\pi }{2}\right]
\end{equation*}
is also conducted. This time, the geometry of the transformation is more complicated. Note that we cannot do $p:[ 0,2\pi ]$ because this will make the mapping (\ref{rvSHvBdfv8WquGXG0ixV}) non-injective, while the model is invertible. Thus the model will have difficulty learning the data near the two lines $p=0$ and $p=2\pi $. The training went for 40,000 epochs, and the loss dropped from $0.3$ to about $10^{-5}$, and the plot is shown in Figure \ref{fig:123203pi2}, and the loss decay is shown in Figure \ref{fig:123203pi2_loss}. The majority of the error comes from $p=3\pi /2$ boundary. This is because the points her are close to the points with $p=0$.


\end{document}